\documentclass[12pt, reqno]{amsart}
\usepackage{amsmath}
\usepackage{amssymb}
 2

\theoremstyle{remark}

\newcommand{\Z}{{\mathbb Z}}
\newcommand{\Q}{{\mathbb Q}}

\newcommand{\C}{{\mathbb C}}

\begin{document}

\title{On the dimension of Chowla-Milnor space }
\author{Tapas Chatterjee}
\address[Tapas Chatterjee]
        {The Institute of Mathematical Sciences, 
         CIT Campus, Taramani, 
         Chennai 600 113 
         India}
\email[Tapas Chatterjee]{tapasc@imsc.res.in}

\maketitle

\begin{abstract}
In a recent work, Gun, Murty and Rath defined the Chowla-Milnor space and proved
a non-trivial lower bound for these spaces. They also  obtained a conditional improvement
of this lower bound and noted that an unconditional  improvement of their lower bound will
lead to irrationality of $\zeta(k)/ \pi^k$ for odd positive integers $k>1$. In this paper,
we give an alternate  proof of their theorem about the conditional lower bound.
\end{abstract}

\section*{Introduction}

\bigskip

For any complex number $s\in\C$, with $\Re(s)>1$, one defines the Riemann zeta function as
\begin{equation*}
 \zeta(s)=\sum_{\substack{n=1}}^\infty \frac{1}{n^s}
\end{equation*}
which has an Euler product 
\begin{equation*}
 \zeta(s)=\prod_{\substack{p}} (1-p^{-s})^{-1}.
\end{equation*}
The Riemann zeta function defines an analytic function in the region 
$\Re(s)>~1$ and can be extended meromorphically to the whole complex plane with a simple pole at $s=1$
having residue 1. Hurwitz generalized the Riemann zeta function by
$\zeta(s,x)$ which is defined as
\begin{equation*}
\zeta(s,x)=\sum_{\substack{n=0}}^\infty \frac{1}{(n+x)^s}
\end{equation*}
where $0<x\leq 1$ and $s\in \C$ with $\Re(s)>1$. He proved that $\zeta(s,x)$ can be extended 
meromorphically to the entire complex plane with a pole at $s=1$. Note 
that for $x=1$, $\zeta(s,1)$ is the classical Riemann zeta function.\\

\noindent
{\bf Definition.} For integers $k >1, q \geq 2$, define the Chowla- Milnor space $V_k(q)$ by
\begin{equation*}
V_k(q) := \Q-{ \rm span ~ of~} \{\zeta(k,a/q):~ 1\leq a <q, ~(a,q)=1\}.
\end{equation*}
 
As described in \cite{GMR}, the conjecture of Chowla and Milnor is the assertion that
the dimension of $V_k(q)$ is equal to $\varphi(q)$, where $\varphi$ is the Euler's phi-function. 
Gun, Murty and Rath \cite{GMR} show that
the dimension of the above spaces is at least $\varphi(q)/2$. They also derived the 
following theorem.

\smallskip
\noindent
{\bf Theorem.} 
Let $k>1$ be an odd integer and $q,r >2$ be two co-prime integers. Then either 
\begin{equation*}
{\rm dim}_\Q V_k(q)\geq \frac{\varphi(q)}{2}+1
\end{equation*}
 or
\begin{equation*}
 {\rm dim}_\Q V_k(r)\geq \frac{\varphi(r)}{2}+1.
\end{equation*}

The proof in \cite{GMR} uses the expansion of Bernoulli polynomials. In this note, we give an
alternate proof of the theorem by an explicit evaluation of 
co-tangent derivatives.

\medskip

\section*{\bf Proof of the Theorem}

\medskip
The following lemma 1 due to Okada~\cite{OK} about the linear independence of co-tangent values 
at rational arguments plays a significant 
role in proving the theorem.

\smallskip
\noindent
{\bf Lemma 1.} Let $k$ and $q$ be positive integers with $k\ge 1$ and $q>2$. Let T be a set of $\varphi(q)/2$
representations mod $q$ such that the union $T\cup(-T)$ constitutes a complete set 
of co-prime residue classes mod $q$.
Then the set of real numbers
\begin{equation*}
 \frac{d^{k-1}}{dz^{k-1}}\cot(\pi z)|_{z=a/q}, ~ a\in T
\end{equation*}
\textit{are linearly independent over $\Q$}.
\smallskip

We first have the following lemma.

\smallskip
\noindent
\textbf{Lemma 2.} \textit{For an integer $k \ge 1$},
\begin{equation*} 
D^{k-1}(\pi\cot\pi z)=\pi^k \times \Z {\rm ~linear ~ combination ~of} ~(\csc\pi z)^{2l}(\cot\pi z)^{k-2l} ,
\end{equation*}
for some non-negative integer $l$. Here $D^{k-1}=\frac{d^{k-1}}{dz^{k-1}}$.

\bigskip
\noindent
\textbf{Proof.}  We will prove this by induction on $k$. For $k=1$, we have\\
$D^{k-1}(\pi\cot(\pi z))=\pi\cot(\pi z)$. Assume that the statement is true for $k-1$, i.e. 
$$
D^{k-2}(\pi\cot(\pi z))=\pi^{k-1}\sum a_i(\csc\pi z)^{2l_i}(\cot\pi z)^{(k-1)-2l_i} 
$$
where $a_i$'s are integers. \\
Differentiating both sides with respect to $z$  we get, 
\begin{eqnarray*}
D^{k-1}(\pi\cot\pi z) &=& \pi^{k}\sum\left[ b_i (\csc\pi z)^{2l_i}
(\cot\pi z)^{k-2l_i}\right. \\
&& \phantom{mmm}+~ \left. c_i(\csc\pi z)^{2l_i+2} (\cot\pi z)^{k-(2l_i+2)}\right],
\end{eqnarray*}
where $b_i, c_i$'s are integers. This completes the proof of lemma 2.

\smallskip
\noindent
{\bf Lemma 3.} \textit{For an integer $k \ge 2$},
\begin{equation*}
\zeta(k,a/q)+(-1)^k\zeta(k,1-a/q)=\frac{(-1)^{k-1}}{(k-1)!}D^{k-1}(\pi\cot\pi z)|_{z=a/q}. 
\end{equation*}

\smallskip
\noindent
{\bf Proof.}
\begin{eqnarray*}
{\rm L.H.S.} &=&  \zeta(k,a/q)+(-1)^k\zeta(k,1-a/q) \\
&=& \underset{n\geq 0}{\sum^\infty} \frac{1}{(n+a/q)^k} ~+ (-1)^k ~\underset{n\geq 0}{\sum^\infty} \frac{1}{(n+1-a/q)^k}\\
&=& \underset{n\geq 0}{\sum^\infty} \frac{1}{(n+a/q)^k} ~+ (-1)^k ~\underset{n=1}{\sum^\infty} \frac{1}{(n-a/q)^k}\\
&=& \underset{n\geq 0}{\sum^\infty} \frac{1}{(n+a/q)^k} ~+ (-1)^{2k} ~\underset{n=1}{\sum^\infty} \frac{1}{(-n+a/q)^k} \\
&=& \underset{n\in\Z}{\sum}\frac{1}{(n+a/q)^k}.
\end{eqnarray*}
Again we know that for $z\notin {\Z}$, 
\begin{eqnarray*}
\pi\cot\pi z &=& \underset{n\in\Z}{\sum}\frac{1}{z+n}. 
\end{eqnarray*}
This implies that
\begin{eqnarray*}
D^{k-1}(\pi\cot\pi z) &=&  (-1)^{k-1} (k-1)!\underset{n\in\Z}{\sum}\frac{1}{(z+n)^k}. 
\end{eqnarray*}
So,
\begin{eqnarray*}
\frac{(-1)^{k-1}}{(k-1)!}D^{k-1}(\pi\cot\pi z)|_{z=a/q} &=&  \underset{n\in\Z}{\sum}\frac{1}{(n+a/q)^k},
\end{eqnarray*}
which completes the proof of lemma 3.

\smallskip
\noindent
Finally, we have lemma 4, whose proof is standard.

\smallskip
\noindent
{\bf Lemma 4.} Let $\rm P$ be the set of primes. We have
\begin{equation*}
\zeta(k)\underset{p \in {\rm P}, \atop p|q}{\prod}(1-p^{-k})=
q^{-k}\sum_{\substack{a=1\\ (a,q)=1}}^{q-1}\zeta(k,a/q).
\end{equation*}
\smallskip

\bigskip
\noindent
\textbf{Proof of the Theorem.} First note that the space $V_k(q)$ is 
also spanned by the following sets of real numbers:
\begin{equation*}
\{\zeta(k,a/q)+\zeta(k,1-a/q)|(a,q)=1, ~ 1\leq a<q/2\}, 
\end{equation*}
\begin{equation*}
 \{\zeta(k,a/q)-\zeta(k,1-a/q)|(a,q)=1, ~ 1\leq a<q/2\}.
\end{equation*}
Now from the lemma 3, we have the following 
\begin{equation*}
\zeta(k,a/q)+(-1)^k\zeta(k,1-a/q)=\frac{(-1)^{k-1}}{(k-1)!}D^{k-1}(\pi\cot\pi z)|_{z=a/q}.
\end{equation*}
Applying the above lemma 1, we see that
\begin{equation*}
 {\rm dim}_\Q V_k(q)\geq \frac{\varphi(q)}{2}.
\end{equation*} 
Now from lemma 2 and lemma 3 for an odd integer $k$, we have 
\begin{eqnarray*}
&& \frac{\zeta(k,a/q)-\zeta(k,1-a/q)}{(2\pi i)^k}   \\
&& =\frac{i}{2^k}\times \Q ~{\rm linear ~ combinations ~of}
~(\csc\pi a/q)^{2l}(\cot\pi a/q)^{k-2l}.
\end{eqnarray*}
We note that $$i \cot(\pi a/q) = \frac{1+\zeta_q^a}{1-\zeta_q^a} $$
belongs to $\Q(\zeta_q)$ and hence so do the numbers $\csc(\pi a/q)^{2l}$
and $\cot(\pi a/q)^{2l}$. Since $k$ is odd, we have 
\begin{equation}\label{1}
 \frac{\zeta(k,a/q)-\zeta(k,1-a/q)}{(2\pi i)^k}\in \Q(\zeta_q)
\end{equation}

\noindent
Now we go back to the main part of the proof. Let $q$ and $r$ be two co-prime integers. Suppose that 
\begin{equation*}
 {\rm dim}_\Q V_k(q)= \frac{\varphi(q)}{2}.
\end{equation*}
Then the numbers 
\begin{equation*}
 \zeta(k,a/q)-\zeta(k,1-a/q), ~ {\rm where} ~(a,q)=1, ~ 1\leq a<q/2
\end{equation*}
generate $V_k(q)$. Now from lemma 4, we get
\begin{equation*}
\zeta(k)\underset{p|q}{\prod}(1-p^{-k}) ~=~
q^{-k}\sum_{\substack{a=1,\\(a,q)=1}}^{q-1}\zeta(k,a/q) \in V_k(q).
\end{equation*}
and hence 
\begin{equation*}
\zeta(k)=\sum_{\substack{(a,q)=1 \\ 1\leq a<q/2}}\lambda_a[\zeta(k,a/q)-\zeta(k,1-a/q)], ~ \lambda_a\in \Q 
\end{equation*}
so that
\begin{equation*}
\frac{\zeta(k)}{(2\pi i)^k}=\sum_{\substack{(a,q)=1\\1\leq a<q/2}}\frac{\lambda_a[\zeta(k,a/q)-\zeta(k,1-a/q)]}{(2\pi i)^k} 
\end{equation*}

\noindent
Thus by \eqref{1}
\begin{equation*}
\frac{\zeta(k)}{i\pi^k}\in \Q(\zeta_q).
\end{equation*}
Similarly , if 
\begin{equation*}
 {\rm dim}_\Q V_k(r)= \frac{\varphi(r)}{2},
\end{equation*}
then
\begin{equation*}
 \frac{\zeta(k)}{i\pi^k}\in \Q(\zeta_r) 
\end{equation*}
and hence 
\begin{equation*}
\frac{\zeta(k)}{i\pi^k}\in\Q(\zeta_q)\cap \Q(\zeta_r).
\end{equation*}
Since any non-trivial finite extension of $\Q$ is ramified, if $\Q(\zeta_q)\cap \Q(\zeta_r)\neq \Q $ then 
there exists a prime which is ramified in $\Q(\zeta_q)\cap \Q(\zeta_r)$, hence both in 
$\Q(\zeta_q)$ and $\Q(\zeta_r)$. Note a prime which ramifies 
in this intersection must necessarily divide both $q$ and $r$. This is impossible because $(q,r)=1$.
So $\Q(\zeta_q)\cap \Q(\zeta_r)=\Q$.
Hence we arrive at a contradiction as $\frac{\zeta(k)}{\pi^k}$ is a real number. Thus 
$$ 
{\rm dim}_\Q V_k(q)\geq \frac{\varphi(q)}{2}+1
\phantom{m}{\rm or}~\phantom{m} {\rm dim}_\Q V_k(r)\geq \frac{\varphi(r)}{2}+1.
$$
This completes the proof of the theorem.

\bigskip
\noindent
{\bf Acknowledgement.}
I would like to thank Sanoli Gun for helpful discussions.

\end{document}